\documentclass[reqno,12pt,a4paper]{amsart}
\usepackage{graphicx,color,amssymb,latexsym,amsmath,amsfonts}
\usepackage{mathrsfs}
\usepackage{dsfont}
\usepackage{pstricks}
\usepackage{pst-char}
\usepackage{pst-slpe}
\usepackage{pst-blur}
\usepackage{graphicx}

\newtheorem{lemma}{\bf Lemma}%[section]

\newtheorem{remark}{Remark}

\newtheorem{theorem}{\bf Theorem}%[section]

\numberwithin{equation}{section} \theoremstyle{plain}
\theoremstyle{definition}

\DeclareMathOperator*{\argmax}{argmax}

\input cyracc.def

\begin{document}
\title[Cramer theorem]
{Cramer's theorem for nonnegative multivariate point processes
with independent increments}
\author{Klebaner F.}
\address{School of Mathematical Sciences, Building 28M, Monash
University, Clayton Campus, Victoria 3800, Australia.}
\email{fima.klebaner@sci.monash.edu.au}

\author{R. Liptser}
\address{Department of Electrical Engineering Systems,
Tel Aviv University, 69978 Tel Aviv, Israel}
\email{liptser@eng.tau.ac.il}

\date{October 17, 2006}
\maketitle
\begin{abstract}
We consider a continuous time version of Cramer's theorem
with nonnegative summands $ S_t=\frac{1}{t}\sum_{i:\tau_i\le
t}\xi_i, \ t \to\infty, $ where $(\tau_i,\xi_i)_{i\ge 1}$ is a
sequence of random variables such that $tS_t$ is a random process
with independent increments.
\end{abstract}

%\jelclass{C13}

%\subclass{60G35,60G51,62G05,62M20,91B70}

\section{\bf Introduction and main result}
\label{sec-1}

The following version of the Cramer theorem  \cite{Cr} can be
extracted from Dembo and Zeitouni, \cite{DZ}.
\begin{theorem}\label{theo-1}
Let  $(\xi_i)_{i\ge 1}$ be a sequence of
nonnegative identically distributed and independent random
variables with $\xi_1$  admitting  the Laplace transform{\rm :}
$$
\mathscr{L}(\lambda)=\mathsf{E}e^{\lambda\xi_1}, \
\lambda\in(-\infty,\Lambda), \quad \exists \
\Lambda=\inf\{\lambda>0:\mathscr{L}(\lambda)=\infty\}.
$$

Then, the family
$$
S_n=\frac{1}{n}\sum_{i=1}^n\xi_i, \quad n\to\infty
$$
obeys the Large Deviation principle {\rm (}LDP{\rm )} in the metric space
$(\mathbb{R}_+,\varrho)$ {\rm (}for the Euclidean metric $\varrho${\rm )} with the rate
$\frac{1}{n}$ and the rate function
$$
I(u)=
\begin{cases}
\sup\limits_{\lambda\in(-\infty,\Lambda)}\big[\lambda
u-g(\lambda)\big], & u>0
\\
-\log\mathsf{P}(\xi_1=0), & u=0,
\end{cases}
$$
where $g(\lambda)$ is the log moment generation
function,
$$
g(\lambda)=\log\mathsf{E}e^{\lambda \xi_1}, \ \lambda<\Lambda.
$$
\end{theorem}

\medskip
In this paper, we study a ``continuous time version'' of this theorem.
For $t\ge 0$, set
\begin{equation*}
S_t=\frac{1}{t}\sum_{i:\tau_i\le t}\xi_i,
\end{equation*}
where $(\xi_i,\tau_i)_{i\ge 1}$ is a sequence of random pairs,
where $\xi_i$'s and $\tau_i$'s are random variables:
$$
\xi_i\ge 0\quad\text{and}\quad \tau_0=0<\tau_1<\tau_2<\ldots
\tau_i<\ldots
$$
defined on a probability space $(\varOmega,\mathcal{F},\mathsf{P})$.
Let $(\mathscr{G}_n)_{n\ge 0}$ be the filtration with
$\mathscr{G}_0=(\varnothing,\varOmega)$ and
$\mathscr{G}_n:=\sigma\{(\tau_i,\xi_i)_{i\le n}\}$.
Random variables $\xi_i$ are assumed to be identically distributed and
independent of $\mathscr{G}_{i-1}$ with the distribution function
$$
G(x)=\mathsf{P}(\xi_1\le x), \ x\ge 0.
$$
The conditional distribution of $\tau_i$ given $\mathscr{G}_{i-1}$ is exponential:
$$
\mathsf{P}\big(\tau_i\le t|\mathscr{G}_{i-1}\big) =
\big(1-e^{-r(t-\tau_{i-1})}\big), \ t\ge \tau_{i-1},
$$
where $r$ is a positive number. Moreover, we assume that
\begin{equation*}
\mathsf{P}\big(\xi_i\le x,\tau_i\le t
|\mathscr{G}_{i-1}\big)=G(x)
\big(1-e^{-r(t-\tau_{i-1})}\big), \ i\ge 1.
\end{equation*}

The following theorem is an analogue of Theorem \ref{theo-1}.
\begin{theorem}\label{theo-2}
The family
\begin{equation*}
S_t=\frac{1}{t}\sum_{i:\tau_i\le t}\xi_i, \quad t\to\infty
\end{equation*}
obeys the LDP in the metric space
$(\mathbb{R}_+,\varrho)$ with the rate
$\frac{1}{t}$ an the rate function
$$
I(u)=
\begin{cases}
\sup\limits_{\lambda\in(-\infty,\Lambda)}\big[\lambda
u-r\int_0^\infty(e^{\lambda z}-1)dG(z)\big] ,& u>0
\\
r[1-G(0+)], & u=0.
\end{cases}
$$
\end{theorem}

\begin{figure}
  \includegraphics[angle=0,width=3.0in,height=2in]{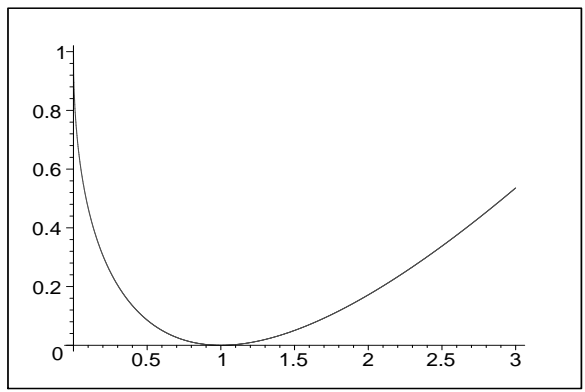}\\
  \caption{The rate function $I^c(u)$ for $r=1$}
  \label{fig-2}
\end{figure}

\medskip
We give two examples illustrating compatibility with Theorems \ref{theo-1} and
\ref{theo-2}.
For both discrete and continuous time cases, let
$$
\mathsf{P}(\xi_1\le x)=1-e^{-x}, \ x\ge 0,
$$
so that, $\xi_1$ has the Laplace transform with $\Lambda=1$. The
log moment generating function is
$$
g(\lambda)=-\log(1-\lambda), \ \lambda<1.
$$
For $G(x)=1-e^{-x}$, $x\ge 0$.
$$
\int_0^\infty r(e^{\lambda z}-1)dG(z)=\frac{r\lambda}{1-\lambda}, \ \lambda<1.
$$
In both cases, rate functions are explicitly computable (see also figures
\eqref{fig-1} and  \eqref{fig-2}),
\begin{equation*}
\begin{aligned}
&
I^d(u)=
  \begin{cases}
    u-1-\log(u), & u>0
    \\
    \infty , & u=0
  \end{cases}
\ \text{(discrete time case)}
\\
&
I^c(u)=
  \begin{cases}
    \big(\sqrt{r}-\sqrt{u}\big)^2 , & u>0
    \\
    r , & u=0.
  \end{cases}
\ \text{(continuous time case)}
\end{aligned}
\end{equation*}

\begin{figure}
\includegraphics[angle=0,width=3.0in,height=2in]{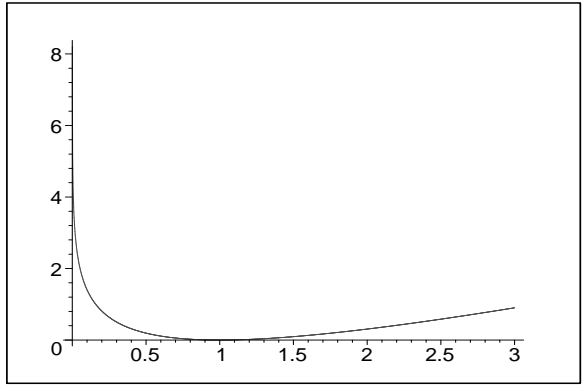}
\label{fig-1}
\caption{The rate function $I^d(u)$}
\label{fig-1}
\end{figure}

\begin{remark}
{\rm Related topics  to Theorem \ref{theo-2} can be found
 e.g. in Georgii and Zessin, \cite{GZ}, serving
a class of marked point random fields. Probably, the proof of
Theorem \ref{theo-2} can be adapted with arguments from proofs in
\cite{GZ} provided that many details not related to our setting
have to be omitted  and other ones concerning to the boundary
effect have to be added.

We prefer to give a complete and direct proof of Theorem
\ref{theo-2}.}
\end{remark}

\section{\bf Counting random measure, its compensator, \\ Laplace transform}
\label{sec-2}

We consider $(\xi_i,\tau_i)_{i\ge 1}$ as a
multivariate (marked) point process (see, e.g. \cite{J},
\cite{JS}) with the counting measure
\begin{equation*}
\mu(dt,dy)=\sum_{i\ge 1}\mathbf{I}_{\{ \tau_i<\infty\} }
\delta_{{\{\tau_i,\xi_i\}}}(t,y)dtdy,
\end{equation*}
where $\delta_{{\{\tau_i,\xi_i\}}}$ is the Dirac delta-function
on $\mathbb{R}_{+}\times\mathbb{R}_+$. Parallel to $(\mathscr{G}_n)_{n\ge 0}$,
we introduce one more filtration $(\mathscr{G}_t)_{t\ge 0}$ related to
$(\xi_i,\tau_i)_{i\ge 1}$:
$$
\mathscr{G}_{t}:=\sigma(\mu([0,t']\times\varGamma):t'
\le t,\varGamma\in \mathscr{B}(\mathbb{R}_+)),
$$
where $\mathscr{B}(\mathbb{R}_+)$ is the Borel $\sigma$-algebra on
$\mathbb{R}_+$, and assume that $\mathscr{G}_0$ is augmented by
$\mathsf{P}$-zero sets from $\mathcal{F}$ (notice that, then,
$(\mathscr{G}_t)_{t\ge 0}$ satisfies the general conditions). With
the help of the counting measure $\mu$, one can present $tS_t$ in
a form of a stochastic integral with respect to $\mu$:
\begin{equation}\label{2.3}
tS_t=\int_0^t\int_{x>0}x\mu(ds,dx).
\end{equation}
Then, the Levy measure $\nu(ds,dx)$, related to $\mu$, is explicitly
computed (see, e.g. Theorem III.1.33, \cite{JS})
\begin{equation}\label{2.2}
\nu(ds,dx)=\sum_{i\ge 1}I_{]\hskip-1.5pt]\tau_{i-1},\tau_{i}]\hskip-1.5pt]}
(t)\frac{dG(x)de^{-r(s-\tau_{i-1})}}{e^{-r(s-\tau_{i-1})}}=rdsdG(x)
\end{equation}
and is not random. It is well known (see, e.g. Corollary to
Theorem 1 in \S 4, Ch.4, \cite{LShMar} that under the
deterministic Levy process $tS_t$ is the random process with
independent increments. We recall a  useful property: for any
nonnegative  and $(\mathscr{G}_t)$-predictable function
$f(\omega,x,t)$,
$$
\mathsf{E}\int_0^t\int_{x>0}f(\omega,x,s)\mu(ds,dx)=
\mathsf{E}\int_0^t\int_{x>0}f(\omega,x,s)\nu(ds,dx)
$$
with ``$\infty=\infty$''.

\begin{lemma}\label{lem-3.1}
{\rm [Laplace transform]}
For any $\lambda<\Lambda$ and $t>0$,
\begin{equation*}
\mathsf{E}e^{\lambda tS_t}=e^{rt\int_{x>0}[e^{\lambda z}-1]dG(z)}.
\end{equation*}
\end{lemma}
\begin{proof}
Though the direct computation of Laplace's transform is permissible, we prefer to
apply the stochastic calculus.
The process
$
U_t=e^{\lambda tS_t}
$
has right-continuous piece-wise constant paths
with jumps
$$
\triangle U_s=(U_s-U_{s-})=U_{s-}\int_{x>0}[e^{\lambda x}-1]\mu(\{s\},dx),
$$
so that, for any $t>0$,
\begin{equation*}
U_t= 1+\int_0^t\int_{x>0}U_{s-}[e^{\lambda x}-1]\mu(ds,dx).
\end{equation*}
Since the function $f(\omega,x,s):=U_{s-}[e^{\lambda x}+1]$ is
nonnegative and predictable, the following equality with
$\lambda<\Lambda$ holds true:
$$
E\int_0^t\int_{x>0}U_{s-}[e^{\lambda x}+1]\mu(ds,dx)=
E\int_0^t\int_{x>0}U_{s-}[e^{\lambda x}+1]\nu(ds,dx)\ (<\infty).
$$
Then, we also have
$$
E\int_0^t\int_{x>0}U_{s-}[e^{\lambda x}-1]\mu(ds,dx)=
E\int_0^t\int_{x>0}U_{s-}[e^{\lambda x}-1]\nu(ds,dx) \ (\in\mathbb{R}).
$$
Since $\nu(ds,dx)=rdsdG(x)$, the later provides
$$
(\mathsf{E}U_t)=1+ \int_0^t\int_{x>0}(EU_s)[e^{\lambda x}-1]dG(x)rds.
$$
This can be written in an equivalent form of differential equation
$$
\frac{d(\mathsf{E}U_t)}{dt}=(\mathsf{E}U_{t})\int_{x>0}[e^{\lambda x}-1]dG(x)r
$$
subject to $(\mathsf{E}U_{0})=1$.

Thus, the desired result holds.
\end{proof}

\section{\bf The proof of Theorem \ref{theo-2}}
\label{sec-3}

We verify the necessary and sufficient conditions for the LDP to hold
(for more details, see Puhalskii, \cite{puh2}):

1) exponential tightness,
$$
\lim_{j\to\infty}\varlimsup_{t\to\infty}\frac{1}{t}\log\mathsf{P}\big(S_t\in
\mathbb{R}_+\setminus\mathcal{K}_j\big)=-\infty,
$$

\hskip .18in where $\mathcal{K}_j$'s are compacts increasing to $\mathbb{R}_+$;

\smallskip
2) local LDP, defining the rate function $I(u)$, $u\in\mathbb{R}_+$
$$
\lim_{\delta\to 0}\lim_{t\to\infty}\frac{1}{t}\log\mathsf{P}\big(|S_t-u|\le\delta\big).
=-I(u).
$$

\subsection{The exponential tightness}
 By choosing
$$
K_j=\{x\in\mathbb{R}_+:x\in[0,j]\}
$$
and applying Chernoff's inequality
with parameter $0.5\Lambda$, we find that
$$
\mathsf{P}(S_t>j)\le e^{-0.5\Lambda j+\log\mathsf{E}e^{0.5\Lambda tS_t}}.
$$
By Lemma \ref{lem-3.1},
$$
\mathsf{E}e^{0.5\Lambda tS_t}=e^{rt\int_{x>0}[e^{0.5\Lambda x}-1]dG(z)}
$$
and, therefore,
$$
\frac{1}{t}\log
\mathsf{P}\big(S_t>j\big)\le -0.5\Lambda j+r\int_{x>0}[e^{0.5\Lambda x}-1]dG(x)
\xrightarrow[j\to\infty]{}-\infty
$$
and 1) is done.

\subsection{The local LDP}

We begin with computation of $I(0)$ and prove
\begin{equation}\label{3.0}
\begin{aligned}
&\varliminf_{\delta\to 0}
\varliminf_{t\to\infty}\frac{1}{t}\log\mathsf{P}(S_t\le\delta)\ge-r[1-G(0+)]
\\
&\varlimsup_{\delta\to 0}
\varlimsup_{t\to\infty}\frac{1}{t}\log\mathsf{P}(S_t\le\delta)\le-r[1-G(0+)]
\end{aligned}
\end{equation}
By \eqref{2.2}, $\{tS_t=0\}=\{\mu((0,t]\times\{x>0\})=0$. Consequently
for any $t>0$,
\begin{equation*}
\mathsf{P}(S_t\le \delta)\ge \mathsf{P}(S_t=0)=\mathsf{P}(tS_t=0)
=\mathsf{P}\big(\mu\big((0,t],\{x>0\}\big)=0\big).
\end{equation*}
The counting process $\pi_t:=\mu\big((0,t],\{x>0\}\big)$ has
independent increments and the rate
$$
\mathsf{E}\mu\big((0,t],\{x>0\}\big)=\nu\big((0,t],\{x>0\}\big)=r[1-G(0+)]t.
$$

 It is a
counting process with the compensator $
\nu\big((0,t],\{x>0\}\big)=r[1-G(0+)]t. $ Therefore, by the
Watanabe theorem, \cite{Wat}, $ \pi_t $ is a Poisson process with
parameter $ r[1-G(0+)]. $ Hence, due to   well known property of
the Poisson process
$$
\mathsf{P}(\pi_t=0)=e^{-tr[1-G(0+)]}.
$$
We find that
$$
\frac{1}{t}\log\mathsf{P}(S_t\le\delta)\ge \frac{1}{t}\log\mathsf{P}(\pi_t=0)
=-r[1-G(+)]
$$
implying the lower bond from \eqref{3.0}.

The upper  bound from \eqref{3.0}
is derived with the help of Laplace's transform  with
$0<\lambda<\Lambda$. To this end, we use identity
$$
1=\mathsf{E}\exp\Big(\lambda tS_t-tr\int_{x>0}[e^{\lambda x}-1]dG(x)\Big)
$$
implying the inequality
\begin{equation*}
1\ge \mathsf{E}I_{\{S_t\le\delta\}}\exp\Big(t\Big[\lambda\delta-r\int_{x>0}
[e^{\lambda x}-1]dG(x)\Big]\Big)
\end{equation*}
being equivalent to
$$
\dfrac{1}{t}\log\mathsf{P}(S_t\le\delta)\le -\lambda\delta+
r\int_{x>0}[e^{-\lambda x}-1]dG(x).
$$
Now, passing $t\to \infty$, we obtain the following upper bound
depending on $\delta$ and $\lambda$:
$$
\varlimsup_{t\to\infty}
\frac{1}{t}\log\mathsf{P}(S_t\le\delta)\le -\lambda\delta+
r\int_{x>0}[e^{\lambda x}-1]dG(x).
$$
Now, passing $\delta\to 0$ and $\lambda$ to $-\infty$  we find that
$$
\varlimsup_{\delta\to 0}\varlimsup_{t\to\infty}
\frac{1}{t}\log\mathsf{P}(S_t\le\delta)\le
-r\int_{x>0}dG(x)=-r[1-G(0+)].
$$

\medskip
We continue the proof by checking the formula for $I(u)$ when $u>0$,
i.e.
$$
\begin{aligned}
&\varlimsup_{\delta\to 0}
\varlimsup_{t\to\infty}\frac{1}{t}\log\mathsf{P}(|S_t-u|\le\delta)\le-I(u)
\\
&\varliminf_{\delta\to 0}
\varliminf_{t\to\infty}\frac{1}{t}\log\mathsf{P}(|S_t-u|\le\delta)\ge-I(u),
\end{aligned}
$$
with
$$
I(u)=\sup\limits_{\lambda\in(-\infty,\Lambda)}\big[\lambda
u-r\int_0^\infty(e^{\lambda x}-1)dG(x)\big].
$$
The Laplace transform
\begin{equation*}
1=\mathsf{E}\exp\Big(\lambda tS_t-tr\int_{x>0}[e^{\lambda
x}-1]dG(x)\Big), \ \lambda<\Lambda,
\end{equation*}
implies the inequality
\begin{equation*}
1\ge
\mathsf{E}I_{\{|S_t-u|\le\delta\}}
\exp\Big(-t\delta u+\lambda tu-tr\int_{x>0}[e^{\lambda x}-1]dG(x)\Big)
\end{equation*}
prviding the following upper bound depending on $\lambda$:
\begin{gather*}
\varlimsup_{\delta\to 0}\varlimsup_{t\to\infty}\frac{1}{t}\log
\mathsf{P}(|S_t-u|\le\delta)\le- \Big(\lambda u-r\int_{x>0}[e^{\lambda x}-1]dG(x)\Big).
\end{gather*}
A further minimization  of the right hand side of the above
inequality in $\lambda$ over $(-\infty,\Lambda)$ gives the desired
result.

\medskip
The lower bound proof uses a standard approach of changing
``probability measure''. Denote by
$\lambda^*=\argmax_{\lambda<\Lambda} \big(\lambda
u-r\int_{x>0}[e^{\lambda x}-1]dG(x)\big)$. Since $\lambda^*$ solves
the equation (with $u>0$)
\begin{equation}\label{leq}
u-r\int_{x>0}xe^{\lambda x}dG(x)\big)=0,
\end{equation}
$\lambda^*$ is a proper number strictly less than $\Lambda$. Set
\begin{equation}\label{LLL}
\frak{L}_t(\lambda^*)=\exp\Bigg(\lambda^*tS_t-
\int_0^t\int_{x>0}r[e^{\lambda^* x}-1]dG(x)ds\Bigg).
\end{equation}
First of all we notice that the Laplace transform for $tS_t$ with
$\lambda^*$ guarantees $\mathsf{E}\frak{L}_t(\lambda^*)=1$.
Moreover, taking into account \eqref{2.3} and applying the It\^o
formula to $\frak{L}_t(\lambda^*)$ one can  see that
$(\frak{L}_t(\lambda^*),\mathscr{G}_t)_{t\ge 0}$ is a positive local
martingale with paths from the Skorokhod space
$\mathbb{D}_{[0,\infty)}$. Then a measure
$\widetilde{\mathsf{P}}_t$, defined by $
d\widetilde{\mathsf{P}}_t=\frak{L}_t(\lambda^*)d\mathsf{P}_t, $
where $\mathsf{P}_t$ is a restriction of $\mathsf{P}$ on
$\mathscr{G}_t$, is the probability measure. We introduce the
probability space $(\varOmega,\mathcal{F},
\widetilde{\mathsf{P}}_t)$. Since $\frak{L}_t(\lambda^*)>0$,
$\mathsf{P}$-a.s., not only $\widetilde{\mathsf{P}}_t\ll
\mathsf{P}_t$ but also $ \widetilde{\mathsf{P}}_t\ll\mathsf{P}_t $
with $
d\mathsf{P}_t=\frak{L}^{-1}_t(\lambda^*)\widetilde{\mathsf{P}}_t. $
This property and \eqref{LLL} provide a lower bound
$$
\begin{aligned}
\mathsf{P}\big(|S_t-u|\le\delta\big)&=\mathsf{P}_t\big(|S_t-u|\le\delta\big)
=\int_{\{|S_t-u|\le\delta\}}\frak{L}^{-1}_t(\lambda^*)d\widetilde{\mathsf{P}}_t
\\
&\ge e^{-\lambda^*t\delta-tI(u)}\widetilde{\mathsf{P}}_t\big(|S_t-u|\le\delta\big).
\end{aligned}
$$
Therefore, we find that
\[
\varliminf_{t\to\infty}\frac{1}{t}\log\mathsf{P}\big(|S_t-u|\le\delta\big)\ge
-I(u)-\lambda^*\delta+\varliminf_{t\to\infty}\frac{1}{t}\log
\widetilde{\mathsf{P}}_t\big(|S_t-u|\le\delta\big).
\]
Obviously, the desired lower bound to obtain it is left to prove that
\begin{equation*}
\varliminf_{t\to\infty}
\widetilde{\mathsf{P}}_t\big(|S_t-u|\le\delta\big)=1
\end{equation*}
or, equivalently,
\begin{equation}\label{ooo}
\lim_{t\to\infty}\widetilde{\mathsf{P}}_t\big(|S_t-u|>\delta\big)=0.
\end{equation}

Thus the last step of the proof deal with \eqref{ooo}. To this end,
we show that ($\widetilde{\mathsf{E}}_t$ denotes the expectation relative to
$\widetilde{\mathsf{P}}_t$)
\begin{equation*}
\lim_{t\to\infty}\widetilde{\mathsf{E}}_t|S_t-u|^2=0.
\end{equation*}
Since $\mathsf{E}\frak{L}_t(\lambda)=1$, it holds
\begin{gather*}
0=\frac{\partial^2\mathsf{E}\frak{L}_t(\lambda)}{\partial \lambda^2}|_{\lambda=\lambda^*}
\\
=
t^2\mathsf{E}\Big(S_t-r\int_{\{x>0\}}xe^{\lambda^* x}dG(x)\Big)^2\frak{L}_t(\lambda^*)
-t\underbrace{r\int_{\{x>0\}}x^2e^{\lambda^* x}dG(x)}_{=u \ \text{(see \eqref{leq})}}
\\
=t^2\widetilde{\mathsf{E}}\Big(S_t-r\int_{\{x>0\}}xe^{\lambda^* x}dG(x)\Big)^2
-t u.
\end{gather*}
Hence,
$$
\widetilde{\mathsf{E}}\big(S_t-u\big)^2=\frac{1}{t}r\mathsf{E}\int_{\{x>0\}}x^2
e^{\lambda^* x}dG(x)\xrightarrow[t\to\infty]{}0.
$$
\qed

\end{document}